\documentclass[a4paper,12pt]{article}
\usepackage{indentfirst} 
\usepackage{amssymb,amscd,latexsym}
\usepackage[intlimits,sumlimits]{amsmath}
\usepackage{amsfonts}
\usepackage{graphicx}
\usepackage{amsthm}
\usepackage{amsxtra}
\usepackage[dvips]{epsfig}
\usepackage[dvips]{graphicx,color}
\usepackage{fancyheadings}
\usepackage[cp1251]{inputenc}
\usepackage[english]{babel}
\usepackage{colortbl}
\begin{document}
\begin{center}
 \vspace{10pt}
 \textbf{The Crank-Nicholson type compact difference scheme for a loaded time-fractional Hallaire's equation}

 \vspace{3mm}

\textbf{Anatoly Alikhanov$^a$, Murat Beshtokov$^a$, Mani Mehra$^b$}\\
\vspace{3mm}
$^a$Institute of Applied Mathematics and Automation KBSC RAS, Nalchik, 360000, Russia \\
$^b$Indian Institute of Technology Delhi, New Delhi, 110016, India
\vspace{10pt}
\end{center}

\textbf{Abstract.} In this paper, we study loaded
modified diffusion equation (the Hallaire equation
with the fractional derivative with respect to time). The compact finite difference scheme
of Crank-Nicholson type of higher order is developed for approximating
the stated problem on uniform grids. A priori estimates are obtained
in difference and differential interpretations, from which there follow
uniqueness, stability, and convergence of the solution of the difference
problem to solution of the differential problem with the rate
$\mathcal{O}(h^4+\tau^{2-\alpha}.)$ Proposed theoretical calculations
are confirmed by numerical experiments on test problem.

\section{Introduction} It is well known that filtration of
liquids in porous media \cite{Barl1,Dz2}, heat transfer in a
heterogeneous environment \cite{Rub3,Ting4}, moisture
transfer in soil grounds \cite{Hil5}, \cite[p.137]{Chud6}  lead to
modified diffusion equations (i.e. the Hallaire
equation \cite{Hil5}). For example, the movement of water in
capillary-porous media, to which the soil belongs, can occur under
the influence of the most diverse driving forces. Based on the
analysis of the diffusion mechanism in porous array when the
occurrence of moisture flows under capillary pressure gradient a
nonlinear equation is obtained \cite[p.136]{Hil5}
$$
W_{t}=(D(W)W_{x})_{x},
$$
where $W$ is humidity in fractions of units, $x$ is depth, $t$
is time and $D(W)$ represents diffusivity.

Diffusion model assuming that if at the initial moment a non-uniform
humidity is given, then a flow of moisture from the more wet layer
into less wet layers occurs, often is not justified. To explain when
and under what conditions the movement of moisture forward and
backward takes place, you can use the modified diffusion equations
(Hallaire equation) \cite{Hil5}:
\begin{equation}
W_{t}=(D(W)W_{x}+AW_{xt})_{x}, \label{ur1}
\end{equation}
where $A$ is a positive constant. An equation of the form (\ref{ur1}) is often called pseudoparabolic.
Various boundary value problems for pseudoparabolic equations were
studied in literature  \cite{Col7, Col8, Run9, Vod10, Besh11, Besh12,
Kozh13, Shkh14, Col15, Show16}. The numerical solutions of the loaded
differential equations are discussed by numerous authors \cite{Shkh17, Shkh18,
Shkh19, Shkh20, Al21, Besh22, Besh23}.

Differential equations of fractional order attract more and more
attention of scientists due to the fact that equations of this type
can describe many physical and chemical processes, biological
environments and systems, which are well interpreted as fractals (i.e soil which is most porous). Overview
of the basic theory of fractional differentiation, fractional order
differential equations, methods for their solution and applications
can be found in \cite{Mod24, Pod25, Nukh26, Uchl27, Samk28}.

Consider the modified diffusion equation (Hallaire equation) with fractional
time derivative
\begin{equation}\label{ur2}
\partial_{0t}^{\alpha}u(x,t)=\mathcal{L}u(x,t)+\mu\frac{\partial}{\partial t}\mathcal{L}u+\mathcal{I} u+f(x,t) ,\,\, 0<x<l,\,\,
0<t\leq T,
\end{equation}
\begin{equation}
u(0,t)=0,\quad u(l,t)=0,\quad 0\leq t\leq T, \quad
u(x,0)=u_0(x),\quad 0\leq x\leq l,\label{ur3}
\end{equation}
where
\begin{equation}
\partial_{0t}^{\alpha}u(x,t)=\frac{1}{\Gamma(1-\alpha)}\int\limits_{0}^{t}\frac{\partial
    u(x,\eta)}{\partial\eta}(t-\eta)^{-\alpha}d\eta, \quad 0<\alpha<1,
\label{ur4}
\end{equation}
is Caputo fractional derivative of order $\alpha$ and 
\begin{equation}\mathcal{L}u(x,t)=
\frac{\partial^2 u}{\partial x^2}, \quad \mathcal{I}u=\sum_{k=1}^m q_k(x,t)u(x_k,t),
\label{ur5}
\end{equation}
$\mu>0, \quad |q(x,t)|\leq c_1$, \quad $0<x_1<x_2<...<x_m<l$.

The most common approximation to the fractional derivative
(\ref{ur4}) is the so-called $L1$ method \cite{Old29, Sun30} which
is defined as follows
\begin{equation}
\partial_{0t_{j+1}}^{\alpha}u(x,t)=\frac{1}{\Gamma(1-\alpha)}\sum\limits_{s=0}^{j}\frac{u(x,t_{s+1})-u(x,t_{s})}{t_{s+1}-t_{s}}\int\limits_{t_{s}}^{t_{s+1}}
\frac{d\eta}{(t_{j+1}-\eta)^{\alpha}}+r^{j+1},
\label{ur6}
\end{equation}
where $0=t_0<t_1<\ldots<t_{j+1}$, and $r^{j+1}$ is approximation
error. In the case of a uniform grid $\tau=t_{s+1}-t_s$, for all
$s=0,1,\ldots, j+1$ and it was proved that
$r^{j+1}=\mathcal{O}(\tau^{2-\alpha})$ \cite{Sun31,Lui32, Al33}. In \cite{Shkh34, Chen35, Al36}  using the $L1$ method, difference
schemes for fractional-order Caputo diffusion equations were
constructed, and their stability and convergence were proved.

In contrast to the classical case, due to the non-local properties of the fractional differentiation operator, the algorithms for solving fractional-order equations are rather laborious even in the one-dimensional case. In the transition to two-dimensional and three-dimensional problems, their complexity increases significantly. In this regard, the construction of stable difference schemes of high order of accuracy is a very important task.

In \cite{Col37, Al38, Lel39, Mehra40, Gui41, Sun42, Sun43, Sun44,
Du45, Ghen46, Gao47} compact difference schemes for increasing the
order of the error of approximation in spatial variables were
construted and investigated.

Approximation of the fractional derivative (\ref{ur4}) using the $L1$ method  gives the error of $\mathcal{O}(\tau^{2-\alpha})$, however
a two-point approximation in the time variable of the term $\mu\dfrac{\partial}{\partial t}\mathcal{L}u$ of the right side of equation (\ref{ur2}) at the grid points $t=t_{j+1}$ gives the error $\mathcal{O}(\tau)$. In this regard, in this paper, we propose to approximate the fractional derivative (\ref{ur4}) on the half layers $t=t_{j+\frac{1}{2}}, j=0,1,2,...$, to keep the order in time $\mathcal{O}(\tau^{2-\alpha})$. For the numerical solution of the problem, some efficient compact difference schemes of Crank-Nicholson type (having high order accuracy) are used.

Difference methods for solving boundary value problems for
pseudoparabolic fractional order equations have been studied in
\cite{Besh48}, \cite{Besh49}.

In this paper, we construct a difference analog of the fractional Caputo derivative on the half-layers, which ensures the order of approximation in the time variable no worse than $ \mathcal {O} (\tau ^ {2- \alpha}) $. To increase the order of approximation with respect to the spatial variable, compact difference schemes are used. Thus constructed difference scheme
has order of accuracy $\mathcal{O}(h^4+\tau^{2-\alpha})$. The method of energy inequalities is used to obtain a priori estimates for both differential and difference problems. Numerical experiments to
verify the accuracy and efficiency of the proposed solution algorithm have been carried out.

\section{A Priori Estimate in the Differential Form}

\textbf{Lemma 1. }\cite{Al50} \textit{For any absolutely continuous
on $[0,T]$ function $\upsilon(t)$
 the following inequality holds true
 $$
 \upsilon(t)\partial_{0t}^\alpha \upsilon(t) \geq \frac{1}{2} \partial_{0t}^\alpha \upsilon^2(t), \ \ \ 0<\alpha<1.
 $$
}
\textbf{Theorem 1.} \textit{If $|q_k(x,t)|<c_1$, $k=1,2,...,m$, then for the solution
to problem (\ref{ur2})-(\ref{ur3}) it holds true that
\begin{equation}
D_{0t}^{\alpha-1}\|u\|_0^2+\|u_x\|_0^2 d\tau \leq M(T) \left(
\int_0^t
\|f\|_0^2d\tau+\|u'_0\|_0^2+\|u_0\|_0^2\right),
\label{ur7}
\end{equation}
where $D_{0t}^{\alpha-1}u=\frac{1}{\Gamma(1-\alpha)}\int\limits_0^t
{(t-\eta)^{-\alpha}}{u(x,\eta) d\eta} $ is Riemann-Liouville fractional integral of order $1-\alpha,$ $M(T)>0 $ is a constant depending only on the input data of problem
(\ref{ur2})-(\ref{ur3}).}

\textbf{Proof.} By the method of energy inequalities we find
a priori estimate for solving problem (\ref{ur2})-(\ref{ur3}). To do this, scalar-multiply
equation (\ref{ur2}) on $u:$
\begin{equation}
\left(u,\partial^{\alpha}_{0t}u\right) =\left(u,\mathcal{L}u\right)
+\mu\left(u,\frac{\partial}{\partial t}\mathcal{L}u\right)+\left(u,\mathcal{I}
u\right)+\left(u,f\right),
\label{ur8}
\end{equation}

$(u,v)=\int\limits^l_0uvdx,$ $(u,u)=\int\limits^l_0u^2dx. $

We transform each term in (\ref{ur7}) subject to Lemma 1 and conditions
(\ref{ur3}):
\begin{equation}
\left(u,\partial^{\alpha}_{0t}u\right)\geq
\frac{1}{2}\partial^{\alpha}_{0t}\|u\|_0^2,
\label{ur9}
\end{equation}
\begin{equation}
\left(u,\mathcal{L}u\right)=\left(u_{xx},u\right)=-\|u_x\|_0^2,
\label{ur10}
\end{equation}
\begin{equation}
\mu\left(u,\frac{\partial}{\partial t}\mathcal{L}u\right)=\mu\left(
u_{xxt},u\right)=-\mu\left(
u_x,u_{xt}\right)=-\frac{\mu}{2}\frac{d}{dt}\|u_x\|_0^2,
\label{ur11}
\end{equation}
$$
\left(u,Iu\right)=\left(u,\sum^m_{k=1}q_ku(x_k,t)\right)=\int_0^l
u(x,t)\sum^m_{k=1}q_k u(x_k,t) dx
$$
\begin{equation}
=\sum^m_{k=1}u(x_k,t)\int_0^l q_k(x,t)u(x,t)dx \leq mc_1 l \|u\|_C^2
\leq \frac{mc_1 l^3}{2} \|u_x\|_0^2,
\label{ur12}
\end{equation}
\begin{equation}
\left(u,f\right)\leq \frac{1}{2} \|f\|_0^2 +\frac{1}{2} \|u\|_0^2
\leq \frac{1}{2} \|f\|_0^2 +\frac{l^2
}{4} \|u_x\|_0^2.
\label{ur13}
\end{equation}
Given the transformations obtained, from (\ref{ur7}) we find
\begin{equation}
\partial^{\alpha}_{0t}\|u\|_0^2+\mu
\frac{d}{dt} \|u_x\|_0^2 \leq (mc_1l^3+l^2/2) \|u_x\|_0^2 +
\|f\|_0^2.
\label{ur14}
\end{equation}

Integrate (\ref{ur14}) from $0$ to $t$ with respect to $\tau$
and applying the Gronwall Lemma we find inequality (\ref{ur7}).

From a priori estimates (\ref{ur7}) follow the uniqueness and stability of
the solutions of problem (\ref{ur2})-(\ref{ur3}) on the input data.

\section{Difference analog of the fractional derivative}

On the uniform grid  $\bar{\omega}_\tau =\{t_j=j\tau,
j=0,\ldots,M, T=\tau M\}$ we construct a difference analog of the fractional
derivative of the function $u(t)\in C^2[0,T]$ at fixed points
$t=t_{j+\frac{1}{2}}, \ j = {0,1,....,M-1}.$
$$
\partial^{\alpha}_{0t_{j+\frac{1}{2}}}u(t) =
\frac{1}{\Gamma(1-\alpha)}\int_0^{t_{j+\frac{1}{2}}}
\frac{u'(\eta)}{(t_{j+\frac{1}{2}}-\eta)^\alpha}d\eta
=\frac{1}{\Gamma(1-\alpha)}\sum^j_{s=1} \int_{t_{s-1}}^{t_{s}}
\frac{u'(\eta)}{(t_{j+\frac{1}{2}}-\eta)^\alpha}d\eta
$$
$$
+\frac{1}{\Gamma(1-\alpha)}\int_{t_{j}}^{t_{j+\frac{1}{2}}}
\frac{u'(\eta)}{(t_{j+\frac{1}{2}}-\eta)^\alpha}d\eta.
$$

On the interval $[t_{s-1},t_s] (1\leq s\leq j)$ we write the linear
interpolation $\Pi_{1,s}u(t)$, using two points $(t_{s-1},
u(t_{s-1})), (t_s, u(t_s)),$ then we get
$$
\Pi_{1,s}u(t) = u({t_{s-1}})\frac{t_s-t}{\tau}+
u({t_{s}})\frac{t-t_{s-1}}{\tau}, \quad (\Pi_{1,s}u(t))' =
u_{\bar{t},s}, \quad u(t)-\Pi_{1,s}u(t) =
\frac{u''(\bar{\xi}_s)}{2}(t-t_{s-1})(t-t_{s}),
$$
where $t \in [t_{s-1},t_s], \ \bar{\xi}_s\in(t_{s-1},t_{s}), \
u_{\bar{t},s} = \frac{u(t_{s})-u(t_{s-1})}{\tau}, u_{t,s} = \frac{u(t_{s+1})-u(t_{s})}{\tau}.$

Considering the latter, we find
$$
\partial^{\alpha}_{0t_{j+\frac{1}{2}}}u(t) \approx\frac{1}{\Gamma(1-\alpha)}\sum^j_{s=1} \int_{t_{s-1}}^{t_{s}}
\frac{\Big(\Pi_{1,s}u(\eta)\Big)'}{(t_{j+\frac{1}{2}}-\eta)^\alpha}d\eta+
\frac{u_{t,j}}{\Gamma(1-\alpha)}\int_{t_{j}}^{t_{j+\frac{1}{2}}}
\frac{1}{(t_{j+\frac{1}{2}}-\eta)^\alpha}d\eta
$$
$$
=\frac{1}{\Gamma(2-\alpha)}\sum^j_{s=1}\Big(t^{1-\alpha}_{j-s+\frac{3}{2}}-t^{1-\alpha}_
{j-s+\frac{1}{2}}\Big)u_{\bar t,s}+\frac{\tau^{1-\alpha}}{2^{1-\alpha}\Gamma(2-\alpha)}u_{t,j}
=\frac{\tau^{1-\alpha}}{\Gamma(2-\alpha)}\sum_{s=0}^j
c^{(\alpha)}_{j-s} u_t ^s= \Delta_{0t_{j+\frac{1}{2}}}^\alpha u,
$$
where
$$
 c_{j}^{(\alpha)} =
 \begin{cases}
   \frac{1}{2^{1-\alpha}}, \ j=0,\\
   (j+\frac{1}{2})^{1-\alpha}-(j-\frac{1}{2})^{1-\alpha},\
   j\geq 1.
 \end{cases}
$$

\textbf{Lemma 2.} \textit{For all $\alpha\in(0,1)$ and $u(t)\in
\mathcal{C}^2[0,t_{j+1}]$
$$
|\partial_{0t_{j+\frac{1}{2}}}^{\alpha}u-\Delta_{0t_{j+\frac{1}{2}}}^\alpha
u|\leq \frac{{2}^\alpha M_2^{j+1}}{4\Gamma(2-\alpha)}\left(
\frac{1-\alpha}{2}+1\right)\tau^{2-\alpha},
$$
where $M_2^{j+1} = \max\limits_{0<t<t_{j+1}} |u''(t)|$.}

\textbf{Proof.} Let the following holds:
$\partial_{0t_{j+\frac{1}{2}}}^{\alpha}u-\Delta_{0t_{j+\frac{1}{2}}}^\alpha
u=R_{1}^{j}+R_{j}^{j+\frac{1}{2}}$, where
$$
R_{1}^{j}=\frac{1}{\Gamma(1-\alpha)}\sum\limits_{s=1}^{j}\int\limits_{t_{s-1}}^{t_{s}}\frac{u'(\eta)d\eta}{(t_{j+\frac{1}{2}}-\eta)^{\alpha}}-
\frac{1}{\Gamma(1-\alpha)}\sum\limits_{s=1}^{j}\int\limits_{t_{s-1}}^{t_{s}}\frac{\left({\Pi}_{1,s}u(\eta)\right)'d\eta}{(t_{j+\frac{1}{2}}-\eta)^{\alpha}}
$$
$$
=\frac{1}{\Gamma(1-\alpha)}\sum\limits_{s=1}^{j}\int\limits_{t_{s-1}}^{t_{s}}\left(u(\eta)-{\Pi}_{1,s}u(\eta)\right)'{(t_{j+\frac{1}{2}}-\eta)^{-\alpha}}d\eta
$$
$$
=-\frac{\alpha}{\Gamma(1-\alpha)}\sum\limits_{s=1}^{j}\int\limits_{t_{s-1}}^{t_{s}}\left(u(\eta)-{\Pi}_{1,s}u(\eta)\right){(t_{j+\frac{1}{2}}-\eta)^{-\alpha-1}}d\eta
$$
$$
=-\frac{\alpha}{2\Gamma(1-\alpha)}\sum\limits_{s=1}^{j}\int\limits_{t_{s-1}}^{t_{s}}u''(\bar\xi_s)(\eta-t_{s-1})(\eta-t_{s}){(t_{j+\frac{1}{2}}-\eta)^{-\alpha-1}}d\eta,
$$

$$
R_{j}^{j+\frac{1}{2}}=\frac{1}{\Gamma(1-\alpha)}\int\limits_{t_{j}}^{t_{j+\frac{1}{2}}}\frac{u'(\eta)d\eta}{(t_{j+\frac{1}{2}}-\eta)^{\alpha}}-
\frac{u_{t,j}}{\Gamma(1-\alpha)}\int\limits_{t_{j}}^{t_{j+\frac{1}{2}}}\frac{d\eta}{(t_{j+\frac{1}{2}}-\eta)^{\alpha}}
$$
$$
=\frac{1}{\Gamma(1-\alpha)}\int\limits_{t_{j}}^{t_{j+\frac{1}{2}}}\frac{(u'(\eta)-u_{t,j})d\eta}{(t_{j+\frac{1}{2}}-\eta)^{\alpha}}
$$
$$
=\frac{1}{\Gamma(1-\alpha)}\int\limits_{t_{j}}^{t_{j+\frac{1}{2}}}\frac{\left(u''(\xi_{j+1}^{(1)})\frac{(t_j-\eta)^2}{2\tau}-u''(\xi_{j+1}^{(2)})\frac{(t_{j+1}-\eta)^2}{2\tau}\right)d\eta}{(t_{j+\frac{1}{2}}-\eta)^{\alpha}},
$$
where $\bar\xi_s \in (t_{s-1}, t_{s})$, $s=1, \ldots, j$ and
$\xi_{j+1}^{(1)}\in (t_{j}, \eta)$, $\xi_{j+1}^{(2)}\in (\eta,
t_{j+1})$.

We estimate the errors $R_{1}^{j}$ and $R_{j}^{j+\frac{1}{2}}$:
$$
|R_{1}^{j}|\leq\frac{\alpha
M_2^{j+1}}{2\Gamma(1-\alpha)}\sum\limits_{s=1}^{j}\int\limits_{t_{s-1}}^{t_{s}}(\eta-t_{s-1})(t_{s}-\eta){(t_{j+\frac{1}{2}}-\eta)^{-\alpha-1}}d\eta
$$
$$
\leq\frac{\alpha
M_2^{j+1}\tau^2}{8\Gamma(1-\alpha)}\sum\limits_{s=1}^{j}\int\limits_{t_{s-1}}^{t_{s}}{(t_{j+\frac{1}{2}}-\eta)^{-\alpha-1}}d\eta=
\frac{\alpha
M_2^{j+1}\tau^2}{8\Gamma(1-\alpha)}\int\limits_{0}^{t_{j}}{(t_{j+\frac{1}{2}}-\eta)^{-\alpha-1}}d\eta
$$
$$
=\frac{
M_2^{j+1}\tau^2}{8\Gamma(1-\alpha)}\left(\frac{{2}^\alpha}{\tau^\alpha}-\frac{1}{(j+\frac{1}{2})^\alpha\tau^\alpha}\right)
\leq\frac{ {2}^\alpha
M_2^{j+1}}{8\Gamma(1-\alpha)}\tau^{2-\alpha},
$$

$$
|R_{j}^{j+\frac{1}{2}}|\leq\frac{M_2^{j+1}\tau}{2\Gamma(1-\alpha)}
\int\limits_{t_{j}}^{t_{j+\frac{1}{2}}}\frac{d\eta}{(t_{j+\frac{1}{2}}-\eta)^{\alpha}}=
\frac{{2}^\alpha M_2^{j+1}}{4\Gamma(2-\alpha)}\tau^{2-\alpha}.
$$
Lemma 2 is proved.

\textbf{Lemma 3.} \cite{Sun44} \textit{ Let $v \in
C^6[x_{i-1},x_{i+1}],$ $\xi_i\in (x_{i-1},x_{i+1}), $ then
$$
\frac{v''(x_{i+1})+10v''(x_i)+v''(x_{i-1})}{12}=\frac{v(x_{i+1})-2v(x_i)+v(x_{i-1})}{h^2}+\frac{h^4}{240}v^{(6)}(\xi_i).
$$
}

\section{Stability and Convergence of Compact Finite-Difference Scheme}

To solve the problem, we apply the compact finite difference method. For this purpose, on
uniform grid $\bar \omega_{h\tau}$ to differential problem (\ref{ur2})-(\ref{ur3})
we put in correspondence the difference scheme of the order of approximation
$\mathcal{O}(h^4+\tau^{2-\alpha}):$
\begin{equation}
\Delta_{0t_{j+\frac{1}{2}}}^\alpha \mathcal{H}_h y_i=Y_{\bar xx,i}+\mu y_{\bar
xxt,i}+\mathcal{H}_h\mathcal{I}\bar{Y}+\mathcal{H}_h\varphi_i^{j+1/2}, \quad
i=1, \ldots, N-1, \quad j=0, \ldots, M-1,
\label{ur15}
\end{equation}

\begin{equation}
y^{j+1}_0=y^{j+1}_N=0, \quad j=0, \ldots, M-1, \quad y_i^0=u_0(x_i),
\quad i=0, \ldots, N, \label{ur16}
\end{equation}
where
$$
\mathcal{H}_h y_i^j=\frac{1}{12}
\left(y_{i+1}^j+10y_i^j+y_{i-1}^j\right)=y^j_i+\frac{h^2}{12}y_{\bar{x}x,i}^j,
\quad i=1, \ldots, N-1, \quad
Y^{j+1}=\frac{1}{2}\left(y^{j+1}+y^j\right),
$$
$$
\bar{y}_{i_k} =
\frac{(x_k-x_{i_k})(x_k-x_{i_k+1})(x_k-x_{i_k+2})}{-6h^3}y_{i_k-1}+\frac{(x_k-x_{i_k-1})(x_k-x_{i_k+1})(x_k-x_{i_k+2})}{2h^3}y_{i_k}
$$
$$
+\frac{(x_k-x_{i_k-1})(x_k-x_{i_k})(x_k-x_{i_k+2})}{-2h^3}y_{i_k+1}+\frac{(x_k-x_{i_k-1})(x_k-x_{i_k})(x_k-x_{i_k+1})}{6h^3}y_{i_k+2},
$$
$$
x_{i_k}\leq x_k \leq x_{i_{k+1}}, \ d^{j+1/2}_i=q(x_i,t_{j+1/2}), \ \varphi^{j+1/2}_i=f(x_i,t_{j+1/2}).
$$
In the sequel we will assume that $h<\min\{x_1,1-x_m\}$.

\textbf{Lemma 4.}\cite[p. 171]{Sam51} \textit{If $f_j $ is a
non-decreasing function, $f_j\geq f_{j-1}$ for all $j=1,2, ...$,
then, from inequality
$$
g_{j+1}=c_0\sum\limits^j_{k=1}\tau g_k+ f_j, \quad j=0,1,2, ...,
\quad g_1 \leq f_0.
$$
an evaluation follows:
$$
g_{j+1}=e^{c_0t_{j+1}}f_j,  \quad j=0,1,2, ...
$$}

\textbf{Lemma 5.} \cite{Al38} \textit{If
$g^{j+1}_j>g^{j+1}_{j-1}>...>g^{j+1}_{0}, \ j=0,1,...,M-1,$ then for
any function $v(t)$, defined on the grid $\bar \omega_{\tau}$, the
following inequality holds
$$
{v^{j+1}}_g\Delta^{\alpha}_{0t_{j+1}}v\geq\frac{1}{2}\
_g\Delta^{\alpha}_{0t_{j+1}}(v^2)+\frac{1}{2g^{j+1}_j}\Big( \
_g\Delta^{\alpha}_{0t_{j+1}}v\Big)^2,
$$
$$
{v^{j}}_g\Delta^{\alpha}_{0t_{j+1}}v\geq\frac{1}{2} \
_g\Delta^{\alpha}_{0t_{j+1}}(v^2)-\frac{1}{2\big(g^{j+1}_{j}+g^{j+1}_{j-1}\big)}\Big(
\ _g\Delta^{\alpha}_{0t_{j+1}}v\Big)^2,
$$
where
$_g\Delta^{\alpha}_{0t_{j+1}}v_i=\sum\limits^{j}_{s=0}(v^{s+1}_i-v^s_i)g^{j+1}_s,
\quad
g^{j+1}_s>0, \quad g^1_{-1}=0.$}

It is obvious that the condition $c^{(\alpha)}_0>c^{(\alpha)}_1>c^{(\alpha)}_2>c^{(\alpha)}_3>...>c^{(\alpha)}_j$ is fulfilled
in the case, when $\alpha>\alpha_0=\log_{3}\dfrac{3}{2}$, since for $\alpha\leq \alpha_0,$ $c^{(\alpha)}_0 \leq c^{(\alpha)}_1$.
In order to use Lemma 5, we transform the resulting difference analog of the fractional derivative to the following form
\begin{equation}
\Delta_{0t_{j+\frac{1}{2}}}^\alpha
y=\bar{\Delta}^{\alpha}_{0t_{j+\frac{1}{2}}}y-\gamma
y_t^{j+1}\tau^{1-\alpha}, \label{ur17}
\end{equation}
where
$$
\bar{\Delta}^{\alpha}_{0t_{j+\frac{1}{2}}}y=\frac{\tau^{1-\alpha}}{\Gamma(2-\alpha)}\sum_{s=0}^j
\bar c _{j-s} ^{(\alpha)} y_t^s,
$$
$$
\begin{cases}
\bar c^{(\alpha)}_0=1,\\
\bar c^{(\alpha)}_{j}=c^{(\alpha)}_j,\ j\geq 1,
\end{cases} \quad \bar c^{(\alpha)}_0>\bar c^{(\alpha)}_1 >\bar c^{(\alpha)}_2>...>\bar c^{(\alpha)}_j,\
j=1,2..., \quad
\gamma=\frac{2^{1-\alpha}-1}{2^{1-\alpha}\Gamma(2-\alpha)}>0.
$$

\textbf{Theorem 2.} \textit{Let conditions $|q_k(x,t)|<c_1$,
$k=1,2,...,m$ be fulfilled, then there is such $\tau_0,$ that if
$\tau\le \tau_0,$ then for the solution of difference problem
(\ref{ur15}), (\ref{ur16}) the following a priori estimate is valid
\begin{equation}
\|y^{j+1}_{\bar x}]|_0^2\leq
M(T)\Big(\sum^{j}_{j\,'=0}\|\mathcal{H}_h\varphi\|_0^2\tau +
\|y^0\|^2_0+\|y^0_{\bar{x}}]|^2_0 \Big), \label{ur18}
\end{equation}
where $M(T)>0 $ is a constant which does not depend on $h$ and $\tau$.}

\textbf {Proof.} Multiply the equation (\ref{ur15}) scalarly by
$y^{j+1}:$
$$
\left(y^{j+1}, \Delta_{0t_{j+\frac{1}{2}}}^\alpha (\mathcal{H}_h
y)\right)=\left(y^{j+1},Y_{\bar xx}\right)+\mu\left(y^{j+1},y_{\bar
xxt}\right)+
$$
\begin{equation}
+\left(y^{j+1},\mathcal{I}\bar{Y}\right)+\left(y^{j+1},\mathcal{H}_h\varphi_i^{j+\frac{1}{2}}\right),
\label{ur19}
\end{equation}

where
$$
\left(y,v\right)=\sum_{i=1}^{N-1}y_iv_ih, \quad
\left(y,y\right)=\|y\|_0^2, \quad \left(\mathcal{H}_hy,y\right)=\|y\|_{\mathcal{H}_h}^2, \quad \frac{2}{3}\|y\|_{0}^2\leq\|y\|_{\mathcal{H}_h}^2\leq\|y\|_{0}^2.
$$
We transform the sums entering into (\ref{ur19}), taking into
account Lemma 5, inequality $\|y\|^2_C \leq \dfrac{l^2}{2}\|y_{\bar x}]|^2,
$ and the conditions (\ref{ur16})
$$
\left(y^{j+1},
\Delta_{0t_{j+\frac{1}{2}}}^\alpha \mathcal{H}_hy
\right)=\left(y^{j+1},\bar
\Delta_{0t_{j+\frac{1}{2}}}^\alpha \mathcal{H}_hy \right) -\gamma\tau^{1-\alpha}\left(y^{j+1},\mathcal{H}_hy^{j+1}_t\right)
$$
$$
=\left(y^{j+1},\bar
\Delta_{0t_{j+\frac{1}{2}}}^\alpha \mathcal{H}_hy \right) -\gamma\tau^{1-\alpha}\left(y^{j+1},y^{j+1}_t\right)-\frac{\gamma\tau^{1-\alpha}h^2}{12}\left(y^{j+1},y^{j+1}_{\bar xxt}\right)
$$
$$
=\left(y^{j+1},\bar
\Delta_{0t_{j+\frac{1}{2}}}^\alpha \mathcal{H}_hy \right) -\frac{\gamma\tau^
    {1-\alpha}}{2}\left(1,( y^2)_t+\tau y_t^2\right)+\frac{\gamma\tau^{1-\alpha}h^2}{12}\left(y^{j+1}_{\bar x
},y^{j+1}_{\bar xt}\right]
$$

$$
=\left(y^{j+1},\bar
\Delta_{0t_{j+\frac{1}{2}}}^\alpha \mathcal{H}_hy \right) -\frac{\gamma\tau^
    {1-\alpha}}{2}\left(1,( y^2)_t+\tau y_t^2\right)+\frac{\gamma\tau^{1-\alpha}h^2}{12}\left(1,( y^2_{\bar x})_t+\tau y_{\bar xt}^2\right]
$$
$$
\geq \frac{1}{2} \bar \Delta_{0t_{j+\frac{1}{2}}}^\alpha
\|y\|_{\mathcal{H}_h}^2 -\frac{\gamma\tau^
    {1-\alpha}}{2}\left((\|y\|^2_{0})_t+\tau \|y_t\|_0^2\right)+\frac{\gamma\tau^{1-\alpha}h^2}{12}( \|y_{\bar x}]|^2_{0})_t
$$
\begin{equation}
\geq \frac{1}{2} \bar \Delta_{0t_{j+\frac{1}{2}}}^\alpha
\|y\|_{\mathcal{H}_h}^2 -\frac{\gamma\tau^
    {1-\alpha}}{2}\left((\|y\|^2_{0})_t+\frac{ l^2\tau}{2} \|y_{\bar xt}]|_0^2\right)+\frac{\gamma\tau^{1-\alpha}h^2}{12}( \|y_{\bar x}]|^2_{0})_t,
\label{ur20}
\end{equation}
$$
\left(y^{j+1},Y_{\bar xx}\right)=-\left(y_{\bar
    x}^{j+1},Y^{j+1}_{\bar x}\right]=-\frac{1}{2}\left(y_{\bar x}^{j+1},y_{\bar x}^{j+1}+y_{\bar
x}^j\right]
$$
\begin{equation}
=-\frac{1}{4} \left(1,(y_{\bar
    x}^{j+1}+y_{\bar x}^j)^2+(y_{\bar x}^{j+1})^2-(y_{\bar
    x}^j)^2\right)=-\frac{1}{4}\|Y^{j+1}_{\bar x}]|_0^2-\frac{\tau}{4}\left(\|y_{\bar
x}]|_0^2\right)_t \label{ur21}
\end{equation}
\begin{equation}
\mu\left(y^{j+1},y_{\bar{x} xt}\right)= -\mu\left( y_{\bar
x}^{j+1},y_{\bar xt}\right]=-\frac{\mu}{2}\left(1,( y^2_{\bar
x})_t+\tau y_{\bar xt}^2\right]= -\frac{\mu}{2}\left((\|y_{\bar
x}]|_0^2)_t+\tau\|y_{\bar x t}]|_0^2\right), \label{ur22}
\end{equation}
$$
\left(\mathcal{I}\bar{Y},y^{j+1}\right) \leq
\sum^m_{k=1}\bar{Y}^{j+1}_{k}\left(d^{j+\frac{1}{2}}_k,y^{j+1}\right)\leq
m c_1 l\|Y^{j+1}\|_C\|y^{j+1}\|_C
$$
\begin{equation}
\leq m c_1
l\left(\|y^{j+1}\|^2_C+\|y^{j+1}\|_C\|y^{j}\|_C\right)\leq\frac{3 m
c_1 l^3}{4}\left( \|y^{j+1}_{\bar x}]|_0^2+\|y_{\bar
x}^{j}]|_0^2\right), \label{ur23}
\end{equation}
\begin{equation}
\left(y^{j+1},\mathcal{H}_h\varphi\right)\leq
\frac{1}{2}\|y^{j+1}]|_0^2+\frac{1}{2}\|\mathcal{H}_h\varphi\|_0^2\leq
\frac{1}{4}\|y_{\bar
x}^{j+1}]|_0^2+\frac{1}{2}\|\mathcal{H}_h\varphi\|_0^2. \label{ur24}
\end{equation}

Taking into account the obtained estimates
(\ref{ur20})$-$(\ref{ur24}), from (\ref{ur19}) we find
$$
\bar \Delta_{0t_{j+\frac{1}{2}}}^\alpha
\|y\|_{\mathcal{H}_h}^2 +\left(\mu+\frac{\tau}{2}+\frac{\gamma\tau^{1-\alpha}h^2}{6}\right)( \|y_{\bar x}]|^2_{0})_t
+\frac{1}{2}\|Y^{j+1}_{\bar x}]|_0^2+\tau\left(\mu- \frac{\gamma l^2}{2}\tau^
{1-\alpha}\right)\|y_{\bar x t}]|_0^2
$$
\begin{equation}
\leq \gamma\tau^
{1-\alpha}(\|y\|^2_{0})_t+\frac{3 m c_1 l^3+1}{2}\left( \|y^{j+1}_{\bar x}]|_0^2+\|y_{\bar
    x}^{j}]|_0^2\right)+\|\mathcal{H}_h\varphi\|_0^2.
\label{ur25}
\end{equation}
Summing (\ref{ur25}) from 0 to $j$ multiply $\tau$, we get
$$
\sum^{j}_{s=0}\bar \Delta_{0t_{j+\frac{1}{2}}}^\alpha
\|y^{s+1}\|_{\mathcal{H}_h}^2\tau+\left(\mu+\frac{\tau}{2}+\frac{\gamma\tau^{1-\alpha}h^2}{6}\right)\|y^{j+1}_{\bar x}]|_0^2+\sum^{j}_{s=0}\|Y^{s+1}_{\bar x}]|_0^2\tau
$$
\begin{equation}
+\tau\left(\mu- \frac{\gamma l^2}{2}\tau^
{1-\alpha}\right)\sum^{j}_{s=0}\|y^{s+1}_{\bar x t}]|_0^2\tau\leq M_{1}\left(\sum^{j}_{s=0}\left( \|y^{s+1}_{\bar x}]|_0^2+\|y_{\bar
    x}^{s}]|_0^2\right)\tau+
\sum^{j}_{s=0}\|\mathcal{H}_h\varphi^{s}\|_0^2\tau +
\|y^0_{\bar{x}}]|^2_0 \right), \label{ur26}
\end{equation}
where $M_1=M_1(\alpha, m,c_1,l,\gamma)$ is a known positive constant.

We estimate the first term on the right-hand side of (\ref{ur26}) as
follows:
$$
\sum^{j}_{s=0}\left( \|y^{s+1}_{\bar x}]|_0^2+\|y_{\bar
    x}^{s}]|_0^2\right)\tau\leq
\sum^{j}_{s=0}\|y^{s+1}_{\bar
x}]|_0^2\tau+\sum^{j}_{s=0}\|y^{s}_{\bar x}]|_0^2\tau\leq
\tau\|y^{j+1}_{\bar x}]|_0^2+2\sum^{j}_{s=1}\|y^{s}_{\bar
x}]|_0^2\tau+\tau\|y^{0}_{\bar x}]|_0^2.
$$
Considering the latter, from (\ref{ur26}) we find
$$
\frac{\tau^{1-\alpha}}{\Gamma(2-\alpha)}\sum^{j}_{s=0} \bar c^{(\alpha)}_{j-s}
\|y^{s+1}\|_{\mathcal{H}_h}^2+\left(\mu+\frac{\tau}{2}+\frac{\gamma\tau^{1-\alpha}h^2}{6}-M_1\tau\right)\|y^{j+1}_{\bar x}]|_0^2+\sum^{j}_{s=0}\|Y^{s+1}_{\bar x}]|_0^2\tau
$$
\begin{equation}
+\tau\left(\mu- \frac{\gamma l^2}{2}\tau^
{1-\alpha}\right)\sum^{j}_{s=0}\|y^{s+1}_{\bar x t}]|_0^2\tau\leq
2M_{1}\sum^{j}_{s=1}\|y^{s}_{\bar x}]|_0^2\tau +
 M_{2}\Big(\sum^{j}_{s=0}\|\mathcal{H}_h\varphi^{s}\|_0^2\tau + \|y^0\|^2_0+\|y^0_{\bar{x}}]|^2_0  \Big), \label{ur27}
\end{equation}
where $M_2 $ is a known positive constant.

Choosing
$\tau\leq\tau_0=\min\left\{\frac{\mu}{2M_{1}},\left(\frac{2\mu}{\gamma
l^2}\right)^{\frac{1}{1-\alpha}}\right\}$ and applying Lemma 4 to
inequality (\ref{ur27}), we obtain a priori estimate (\ref{ur18}).

From a priori estimate (\ref{ur18}), uniqueness and stability
follow, and also the convergence of the solution of the difference
problem to the solution of differential problem with the rate
$\mathcal{O}(h^4+\tau^{2-\alpha})$.

\textbf{Remark.} \textit{Obviously, the results obtained in this paper
also hold in the case when
$$
\mathcal{I}u=\int\limits_0^l q(x,t)u(x,t) dx,
$$
if  $|q(x,t)|<c_1.$
}

In this case, in order to preserve the order of approximation of a compact difference scheme, it suffices to apply the Simpson's rules.

\section{Results of a numerical experiment}

Consider the following test case.
$$
\partial^{\alpha}_{0t}u =u_{xx} +\mu u_{xxt}+\sum^3_{s=1}q_s(x,t)u(x_s,t)+ f,
$$
$$
u(0,t)=0, \quad u(1,t)=0,
$$
$$
u(x,0)=u_0(x),
$$
where
the right side $f(x,t)$ is chosen so that the function $u(x,t)=(t^3+t^{2+\alpha}+1)\sin(3\pi x)$ is the exact solution of the original equations. Wherein
$$
\mu=1, \quad x_1=0.2, \quad x_2=0.5, \quad x_3=0.8,
\quad T=1.
$$
$$
q_1(x,t)=e^{x+t},\quad q_2(x,t)=\sin(x+t),\quad q_3(x,t)=\cos(x+t).
$$

The tables 1 and 2, for different values $\alpha$ = 0.1; 0.5; 0.9 and
reducing the size of the grid, shows the maximum value of the error
($z=y-u$) and the convergence order  (CO) in the norms
$\|\,\cdot\,\|_{C(\bar{\omega}_{h\tau})}$ and $\|\cdot\|_{0}$, where
$\|y\|_{C(\bar{\omega}_{h\tau})}=\max\limits_{(x_i,t_j)\in\bar{\omega}_{h\tau}}|y|$.
The error decreases in accordance with the order
of approximations of $\mathcal{O}(h^4+\tau^{2-\alpha})$.  The convergence order
is determined by the following formula:
CO$=\log_{\frac{h_1}{h_2}}\dfrac{\|z_1\|}{\|z_2\|},$
where $z_i$  is the error corresponding to $h_i$.

\vspace{5mm}

\begin{table}[h!]
    \begin{center}
        \caption{The error and the convergence order in the norms $\|\cdot\|_{0}$
            and $\|\cdot\|_{C(\bar{\omega}_{h\tau})}$ when
            decreasing time-grid size for different values of $\alpha=0.1; 0.5; 0.9$ on $t=1$, $\tau=1/10000.$}
        \label{tab:table1}
\begin{tabular}{|c|c|c|c|c|c|c|c|} 
	\hline
	\, \textbf{$\alpha$}\, & \, \textbf{$h$}  \, & \, \textbf{$\|z\|_{C(\bar{\omega}_{h\tau})}$}  \, & \, \textbf{ CO } \, &\, \textbf{$\max\limits_{0\leq j\leq M}\|z^j\|_{0}$}\, & \,\textbf{ CO} \, &\, \textbf{$\max\limits_{0\leq j\leq M}\|z^j_{\bar x}]|_{0}$} \,& \, \textbf{ CO } \,\\
	\hline
		0.1 & 1/6    & 9.757379e-2 &          & 5.430370e-2   &           & 4.132291e-1   &               \\
			& 1/12   & 8.070053e-3 & 3.5958   & 4.024721e-3   &  3.7541   & 2.723858e-2   &  3.9232                 \\
			& 1/24   & 4.202591e-4 & 4.2632   & 2.144845e-4   &  4.2299   & 1.685552e-3   &  4.0143                    \\			
           \hline
		0.5 & 1/6    & 9.661932e-2 &          & 5.404311e-2   &           & 4.120564e-1   &               \\
			& 1/12   & 7.921449e-3 & 3.6085   & 3.952770e-3   &  3.7731   & 2.708703e-2   &  3.9271                 \\
			& 1/24   & 4.140444e-4 & 4.2579   & 2.120782e-4   &  4.2201   & 1.680156e-3   &  4.0109                    \\			
           \hline
		0.9 & 1/6    & 9.561826e-2 &          & 5.375085e-2   &           & 4.105740e-1   &               \\
			& 1/12   & 7.768228e-3 & 3.6216   & 3.879776e-3   &  3.7922   & 2.691586e-2   &  3.9311                 \\
			& 1/24   & 4.075689e-4 & 4.2524   & 2.095680e-4   &  4.2105   & 1.673201e-3   &  4.0077                    \\			
           \hline
        \end{tabular}
    \end{center}
\end{table}

\begin{table}[h!]
    \begin{center}
        \caption{The error and the convergence order in the norms $\|\cdot\|_{0}$
            and $\|\cdot\|_{C(\bar{\omega}_{h\tau})}$ when
            decreasing space-grid size for different values of $\alpha=0.1; 0.5; 0.9$ on $t=1$, $h=1/1000.$}
        \label{tab:table1}
        \begin{tabular}{|c|c|c|c|c|c|c|c|} 
            \hline
           \, \textbf{$\alpha$}\, & \, \textbf{$\tau$}  \, & \,\textbf{$\|z\|_{C(\bar{\omega}_{h\tau})}$} \, & \, \textbf{ CO } \, &\,\textbf{$\max\limits_{0\leq j\leq M}\|z^j\|_{0}$} \, & \,\textbf{ CO} \, &\, \textbf{$\max\limits_{0\leq j\leq M}\|z^j_{\bar x}]|_{0}$} \,& \, \textbf{ CO } \,\\
            \hline
        0.1 & 1/10   & 7.046545e-3 &          & 4.429421e-3   &           & 4.137926e-2   &               \\
            & 1/20   & 1.764596e-3 & 1.9975   & 1.110134e-3   &  1.9963   & 1.037209e-2   &  1.9962                 \\
            & 1/40   & 4.418103e-4 & 1.9978   & 2.781883e-4   &  1.9966   & 2.599478e-3   &  1.9964                    \\
            & 1/80   & 1.106090e-4 & 1.9979   & 6.970348e-5   &  1.9967   & 6.514120e-4   &  1.9965                  \\
            & 1/160  & 2.768977e-5 & 1.9980   & 1.746323e-5   &  1.9969   & 1.632212e-4   &  1.9967                    \\
            & 1/320  & 6.931562e-6 & 1.9981   & 4.374859e-6   &  1.9970   & 4.089443e-5   &  1.9968                    \\
            & 1/640  & 1.734385e-6 & 1.9987   & 1.095364e-6   &  1.9978   & 1.023998e-5   &  1.9976                    \\
            & 1/1280 & 4.327143e-7 & 2.0029   & 2.732386e-7   &  2.0031   & 2.554293e-6   &  2.0032                    \\
            & 1/2560 & 1.066077e-7 & 2.0211   & 6.675703e-8   &  2.0331   & 6.232795e-7   &  2.0349                    \\
            & 1/5120 & 2.604856e-8 & 2.0330   & 1.650392e-8   &  2.0161   & 1.542674e-7   &  2.0144                    \\
            \hline
        0.5 & 1/10   & 8.224209e-3 &          & 5.239804e-3   &           & 4.904309e-2   &               \\
			& 1/20   & 2.067613e-3 & 1.9919   & 1.319392e-3   &  1.9896   & 1.235161e-2   &  1.9893                 \\
			& 1/40   & 5.204621e-3 & 1.9901   & 3.327375e-4   &  1.9874   & 3.115679e-3   &  1.9870                    \\
			& 1/80   & 1.312926e-4 & 1.9870   & 8.412851e-5   &  1.9837   & 7.879783e-4   &  1.9833                  \\
			& 1/160  & 3.322563e-5 & 1.9824   & 2.135228e-5   &  1.9782   & 2.000616e-4   &  1.9777                    \\
			& 1/320  & 8.445963e-6 & 1.9759   & 5.448625e-6   &  1.9704	  & 5.107308e-5   &  1.9698                    \\
			& 1/640  & 2.160852e-6 & 1.9666   & 1.401153e-6   &  1.9593   & 1.314083e-5   &  1.9585                    \\
			& 1/1280 & 5.571053e-7 & 1.9555   & 3.635032e-7   &  1.9466   & 3.411180e-6   &  1.9457                    \\
			& 1/2560 & 1.457657e-7 & 1.9342   & 9.614532e-8   &  1.9186   & 9.030448e-7   &  1.9174                    \\
			& 1/5120 & 3.790931e-8 & 1.9430   & 2.484624e-8   &  1.9522   & 2.332695e-7   &  1.9527                    \\
			\hline
		0.9 & 1/10   & 9.385612e-3 &          & 6.006600e-3   &           & 5.625169e-2   &                       \\
			& 1/20   & 2.369760e-3 & 1.9857   & 1.519540e-3   & 1.9829    & 1.423377e-2   & 1.9825                 \\
			& 1/40   & 6.028051e-4 & 1.9750   & 3.879506e-4   & 1.9697    & 3.635528e-3   & 1.9690                    \\
			& 1/80   & 1.554828e-4 & 1.9550   & 1.007197e-4   & 1.9455    & 9.445138e-4   & 1.9445                \\
			& 1/160  & 4.109394e-5 & 1.9198   & 2.691906e-5   & 1.9036    & 2.527009e-4   & 1.9021                    \\
			& 1/320  & 1.130943e-5 & 1.8614   & 7.542754e-6   & 1.8355    & 7.089996e-5   & 1.8335                    \\
			& 1/640  & 3.310073e-6 & 1.7726   & 2.265946e-6   & 1.7350    & 2.132457e-5   & 1.7332                    \\
			& 1/1280 & 1.051882e-6 & 1.6539   & 7.436887e-7   & 1.6073    & 7.002065e-6   & 1.6066                    \\
			& 1/2560 & 3.876684e-7 & 1.4401   & 2.694090e-7   & 1.4649    & 2.535332e-6   & 1.4656                    \\
			& 1/5120 & 1.578844e-7 & 1.2959   & 1.079825e-7   & 1.3190    & 1.015036e-6   & 1.3206                    \\
\hline
        \end{tabular}
    \end{center}
\end{table}

\section{Conclusion}
In this paper, a difference analogue of the fractional Caputo derivative is constructed on the half-layers, which most effectively approximates the problem under consideration in a time variable. Compact difference schemes are used to increase the order of approximation in the spatial variable. The difference scheme constructed in this paper
has order of accuracy $\mathcal{O}(h^4+\tau^{2-\alpha})$. The stability and convergence of difference schemes with a speed equal to the order of approximation error are proved. The proposed theoretical calculations
are confirmed by numerical experiments on test problems. All numerical calculations were performed using Julia (version 1.0.1).

\center
\end{document}